\documentclass[12pt]{amsart}

\newtheorem{theorem}{Theorem}[section]
\newtheorem{lemma}[theorem]{Lemma}
\newtheorem{proposition}[theorem]{Proposition}
\newtheorem{corollary}[theorem]{Corollary}

\numberwithin{equation}{section}

\begin{document}

\baselineskip=15pt

\title[Restrictions of tangent bundle]{Semistability 
and restrictions of tangent bundle to curves}

\author[I. Biswas]{Indranil Biswas}

\address{School of Mathematics, Tata Institute of Fundamental
Research, Homi Bhabha Road, Bombay 400005, India}

\email{indranil@math.tifr.res.in}

\subjclass[2000]{14F05, 32L10}

\keywords{Semistability, tangent bundle, restriction}

\date{}

\begin{abstract}

We consider all complex projective manifolds $X$
that satisfy at least one of the following three conditions:
\begin{enumerate}
\item There exists a pair $(C\, ,\varphi)$,
where $C$ is a compact connected Riemann surface and
$$
\varphi\,:\, C\,\longrightarrow\, X
$$
a holomorphic map, such that the pull back $\varphi^*TX$ is
not semistable.

\item The variety $X$ admits an \'etale covering by an
abelian variety.

\item The dimension $\dim X\, \leq \,1$.
\end{enumerate}
We prove that the following classes are among those
that are of the above type.
\begin{itemize}
\item All $X$ with a finite fundamental group.

\item All $X$ such that there is a nonconstant morphism from
${\mathbb C}{\mathbb P}^1$ to $X$.

\item All $X$ such that the canonical
line bundle $K_X$ is either positive or
negative or $c_1(K_X)\, \in\, H^2(X,\, {\mathbb Q})$
vanishes.

\item All $X$ with $\dim_{\mathbb C} X\, =\,2$.
\end{itemize}
\end{abstract}

\maketitle

\section{Introduction}

The tangent bundle of a complex projective manifold equipped
with a polarization is often semistable. For example, if
$X$ is a complex projective manifold such that the canonical
line bundle $K_X$
is ample, then the tangent bundle $TX$ is semistable
with respect to the polarization defined by $K_X$.
More generally, if $X$ admits a K\"ahler--Einstein metric
then $TX$ is semistable.
Let $V$ be a holomorphic vector bundle on a
complex projective manifold $X$ equipped with a very
ample line bundle $\zeta$. If $V$ is semistable, then the
restriction of $V$ to any smooth complete intersection
curve in $X$, obtained
by intersecting hyperplanes from the linear
systems of sufficiently large powers of $\zeta$, remains
semistable (see \cite[Ch.~7]{HL}).

Here we consider all connected complex projective manifolds $X$
with the property that for every pair of the form $(C\, ,
\varphi)$, where $C$ is a compact connected Riemann surface and
$$
\varphi\,:\, C\,\longrightarrow\, X
$$
a holomorphic map, the pull back $\varphi^*TX$ is a semistable
vector bundle over $C$. If $\dim_{\mathbb C}X\, \leq\, 1$, or 
$X$ admits an \'etale covering by an abelian variety, then $X$
satisfies this condition (if $A\, \longrightarrow\,X$ is an
\'etale covering with $A$ an abelian variety, then the pull back
of $\varphi^*TX$ to the fiber product $C\times_X A$ is trivial;
hence $\varphi^*TX$ is semistable). We
conjecture that these are all.

For convenience, let us define that a connected complex
projective manifold $M$ satisfies \textit{Condition C}
if at least one of the following three statements holds:
\begin{enumerate}
\item There exists a pair $(Y\, ,\varphi)$,
where $Y$ is a compact connected Riemann surface and
$\varphi\,:\, Y\,\longrightarrow\, M$
a holomorphic map, such that $\varphi^*TM$ is
not semistable.

\item The variety $M$ admits an \'etale covering by an
abelian variety.

\item The variety $M$ is a curve or a point.
\end{enumerate}

The above conjecture says that all connected complex projective
manifolds satisfy Condition \textit{C}. We prove that
the following classes satisfy Condition \textit{C}:
\begin{itemize}
\item All $M$ with a finite fundamental group (Theorem \ref{th1}).

\item All $M$ such that there is a nonconstant morphism from
${\mathbb C}{\mathbb P}^1$ to $M$ (Proposition \ref{prop2}).

\item All $M$ such that either the canonical
line bundle $K_M$ is ample or $K^{-1}_M$ is
ample or $c_1(K_M)\, \in\, H^2(M,\, {\mathbb Q})$
vanishes (Corollary \ref{cor1}).

\item All $M$ with $\dim_{\mathbb C} M\, =\,2$ (Proposition
\ref{prop3}).
\end{itemize}

\section{Flat connection and fundamental group}

Let $M$ be an irreducible smooth complex projective variety.
The complex dimension of $M$ will be denoted by $d$.
Let ${\mathbb P}(TM)$ denote the projectivized tangent bundle
that parametrizes all lines in the tangent spaces of
$M$. Let $F_{\text{PGL}(d,{\mathbb C})}$ be the holomorphic
principal $\text{PGL}(d,{\mathbb C})$ over $M$ defined
by ${\mathbb P}(TM)$.

We recall that a holomorphic connection on
${\mathbb P}(TM)$ is a holomorphic splitting of the
Atiyah exact sequence for the $\text{PGL}(d,
{\mathbb C})$--bundle $F_{\text{PGL}(d,{\mathbb C})}$
(see \cite[page 188, Definition]{At}). The projective bundle
${\mathbb P}(TM)$ admits a flat holomorphic
connection if and only if it admits
local holomorphic trivializations such that all the
transition functions are locally constant.

\begin{proposition}\label{prop1}
Assume that $M$ has the property that for every pair of the
form $(Y\, ,\varphi)$,
where $Y$ is a compact connected Riemann surface and
$$
\varphi\,:\, Y\,\longrightarrow\, M
$$
a holomorphic map, the pull back $\varphi^*TM$ is a semistable
vector bundle over $Y$. Then the projective bundle
${\mathbb P}(TM)$ admits a flat holomorphic connection.
\end{proposition}

\begin{proof}
Consider the adjoint action of $\text{PGL}(d,{\mathbb C})$
on the Lie algebra $\text{M}(d,{\mathbb C})$ of
$\text{GL}(d,{\mathbb C})$. Let
\begin{equation}\label{e1}
\rho\, :\, \text{PGL}(d,{\mathbb C})\, \longrightarrow\,
\text{GL}(\text{M}(d,{\mathbb C}))\, =:\, G
\end{equation}
be the corresponding homomorphism. This homomorphism $\rho$
is injective.

The Lie algebra of $G$ in Eq. \eqref{e1} will be denoted 
by $\mathfrak g$. The Lie algebra of $\text{PGL}(d,{\mathbb C})$
is the subalgebra of $\text{M}(d,{\mathbb C})$ defined by the
trace zero matrices. We will denote the Lie algebra of
$\text{PGL}(d,{\mathbb C})$ by $\text{M}_0(d,{\mathbb C})$ Let
\begin{equation}\label{e2}
\widehat{\rho}\, :\, \text{M}_0(d,{\mathbb C})\,
\longrightarrow\, {\mathfrak g}
\end{equation}
be the injective homomorphism of Lie algebras associated to
$\rho$ in Eq. \eqref{e1}. The group $\text{PGL}(d,{\mathbb C})$
has the adjoint action $\text{M}_0(d,{\mathbb C})$, and
$\text{PGL}(d,{\mathbb C})$ acts on ${\mathfrak g}$ by
combining the homomorphism $\rho$ and the adjoint
action of $G$ on $\mathfrak g$. The injective
homomorphism $\widehat{\rho}$ in Eq. \eqref{e2} is a
homomorphism of $\text{PGL}(d,{\mathbb C})$--modules. The
group $\text{PGL}(d,{\mathbb C})$ is reductive. Hence any
short exact sequence of $\text{PGL}(d,{\mathbb C})$--modules
splits. Therefore, there is a homomorphism of $\text{PGL}(d,
{\mathbb C})$--modules
\begin{equation}\label{e3}
\eta\, :\, {\mathfrak g}\, \longrightarrow\,
\text{M}_0(d,{\mathbb C})
\end{equation}
such that $\eta\circ\widehat{\rho}\, =\, 
\text{Id}_{\text{M}_0(d,{\mathbb C})}$. Fix such a homomorphism
$\eta$.

Let
\begin{equation}\label{e4}
F_G\, :=\, F_{\text{PGL}(d,{\mathbb C})}(G)
\end{equation}
be the holomorphic principal $G$--bundle over $M$ obtained by
extending the structure group of $F_{\text{PGL}(d,{\mathbb C})}$
using the homomorphism $\rho$ in Eq. \eqref{e1}. We recall that
$F_G$ is a quotient of $F_{\text{PGL}(d,{\mathbb C})}\times G$.
Two points $(z_1\, ,g_1)$ and $(z_2\, ,g_2)$ of
$F_{\text{PGL}(d,{\mathbb C})}\times G$ are identified in
$F_G$ if and only if there is an element $g\, \in\, 
\text{PGL}(d,{\mathbb C})$ such that
$$
(z_2\, ,g_2)\,=\, (z_1g^{-1}\, ,\rho(g)g_1)\, .
$$
We have a holomorphic map
\begin{equation}\label{e0}
\alpha\, :\, F_{\text{PGL}(d,{\mathbb C})}\, \longrightarrow\,
F_G
\end{equation}
that sends any $z\, \in\, F_{\text{PGL}(d,{\mathbb C})}$ to the
element in $F_G$ defined by $(z\, ,e)$. Since $\rho$ is injective,
the map $\alpha$ in Eq. \eqref{e0} is an embedding.

We note that the vector bundle over $M$
associated to $F_G$ for the standard action of
$G\, =\, \text{GL}(\text{M}(d,{\mathbb C}))$ on
$\text{M}(d,{\mathbb C})$ is the endomorphism bundle
$$
{\mathcal E}nd(TM)\, =\, TM\bigotimes \Omega^1_M\, .
$$

Fix a very ample line bundle over $M$ to define semistable
vector bundles on it.

Since $\varphi^*TM$ is a semistable vector bundle over $Y$ for
all pairs $(Y\, ,\varphi)$ as in the statement of the proposition,
we know that the vector bundle ${\mathcal E}nd(TM)$ is semistable
and $c_2({\mathcal E}nd(TM))\, \in\, H^4(M,\, {\mathbb Q})$
vanishes (see \cite[Theorem 1.2]{BB}). Also,
$c_1({\mathcal E}nd(TM))\,=\, 0$ because
${\mathcal E}nd(TM)\,=\, {\mathcal E}nd(TM)^*$. Hence the vector
bundle ${\mathcal E}nd(TM)$ admits a flat holomorphic connection
\cite[page 40, Corollary 3.8]{BG} (we note that since $M$
is a complex projective manifold, the semistable vector bundle
${\mathcal E}nd(TM)$ is pseudostable \cite[page 26, Proposition
2.4]{BG}; hence the above mentioned Corollary 3.8 of
\cite{BG} applies).

Giving a flat holomorphic connection on the vector bundle
${\mathcal E}nd(TM)$ is equivalent to giving a flat holomorphic 
connection on the principal $G$--bundle $F_G$ in Eq. \eqref{e4}.
Fix a flat holomorphic connection $\nabla^G$ on the principal
$G$--bundle $F_G$. Therefore,
$$
\nabla^G\, \in\, \Gamma(F_G,\, \Omega^1_{F_G}\bigotimes
{\mathfrak g})
$$
is a holomorphic one--form on the total space of
$F_G$, with values in $\mathfrak g$,
satisfying the following two conditions:
\begin{itemize}
\item the restriction of $\nabla^G$ to any fiber of the
projection $F_G\, \longrightarrow\, M$ coincides with the
Maurer--Cartan form, and

\item the form $\nabla^G$ is equivariant for the action of
$G$ on $F_G$ and the adjoint action of $G$ on $\mathfrak g$.
\end{itemize}

Consider the $\mathfrak g$--valued holomorphic one--form
$\alpha^*\nabla^G$ on $F_{\text{PGL}(d,{\mathbb C})}$, where
$\alpha$ is the embedding constructed in Eq. \eqref{e0}. The
composition $\eta\circ(\alpha^*\nabla^G)$, where $\eta$ is
the projection in Eq. \eqref{e3}, is a 
$\text{M}_0(d,{\mathbb C})$--valued holomorphic one--form
on $F_{\text{PGL}(d,{\mathbb C})}$.

Since $\eta$ in Eq. \eqref{e3} is a homomorphism of
$\text{PGL}(d,{\mathbb C})$--modules satisfying the condition
that
$$
\eta\circ\widehat{\rho}\, =\,
\text{Id}_{\text{M}_0(d,{\mathbb C})}\, ,
$$
it follows that $\eta\circ(\alpha^*\nabla^G)$ defines a
holomorphic connection on the principal $\text{PGL}(d,{\mathbb 
C})$--bundle $F_{\text{PGL}(d,{\mathbb C})}$.

The curvature
of this holomorphic connection on $F_{\text{PGL}(d,{\mathbb C})}$
defined by $\eta\circ(\alpha^*\nabla^G)$
clearly coincides with $\eta\circ{\mathcal K}(\nabla^G)$, where
${\mathcal K}(\nabla^G)$ is the curvature of the connection
$\nabla^G$ on $F_G$. But $\nabla^G$ is flat. Hence
the holomorphic connection on $F_{\text{PGL}(d,{\mathbb C})}$
defined by $\eta\circ(\alpha^*\nabla^G)$ is flat.
This completes the proof of the proposition.
\end{proof}

But we put it down the following lemma for later reference.

\begin{lemma}\label{lem1}
Take $M$ as in Proposition \ref{prop1}. Then
$$
(d-1)c_1(TM)^2\, =\, 2d{\cdot}c_2(TM)\, .
$$
\end{lemma}

\begin{proof}
{}From \cite[Theorem 1.2]{BB} we have
$c_2({\mathcal E}nd(TM))\, =\, 0$. Now the lemma follows from
the fact that $c_2({\mathcal E}nd(TM))\, =\,
2d{\cdot}c_2(TM) -(d-1)c_1(TM)^2$.
\end{proof}

\begin{theorem}\label{th1}
Let $M$ be a connected complex projective manifold
of complex dimension $d$, with $d\, \geq\, 2$. Assume that
$M$ has the property that for every pair of the
form $(Y\, ,\varphi)$,
where $Y$ is a compact connected Riemann surface and
$$
\varphi\,:\, Y\,\longrightarrow\, M
$$
a holomorphic map, the pull back $\varphi^*TM$ is a semistable
vector bundle over $Y$. Then the cardinality of the
fundamental group of $M$ is infinite.
\end{theorem}

\begin{proof}
Assume that the fundamental group of $M$ is finite. Fix
a universal cover
\begin{equation}\label{e5}
\gamma\, :\, \widetilde{M}\, \longrightarrow\, M
\end{equation}
of $M$. Since the fundamental group of $M$ is finite,
this $\widetilde{M}$ is also a connected complex
projective manifold of complex dimension $d$.

Let $Y$ be a compact connected Riemann surface, and let
$$
\phi\, :\, Y\, \longrightarrow\, \widetilde{M}
$$
be a holomorphic map. Set
$$
\varphi\, :=\, \gamma\circ \phi\, ,
$$
where $\gamma$ is the map in Eq. \eqref{e5}. Since $\gamma$ is
an \'etale covering, we have
$$
\varphi^*TM\, =\, \phi^* T\widetilde{M}\, .
$$
Therefore, using the given condition on $M$ it follows
that the vector bundle $\phi^* T\widetilde{M}$ is semistable.

Hence from Proposition \ref{prop1} we know that the
projective bundle ${\mathbb P}(T\widetilde{M})$ admits a flat
holomorphic connection. On the other hand, $\widetilde{M}$
is simply connected. Hence the projective bundle ${\mathbb
P}(T\widetilde{M})$ is trivial. This
immediately implies that
the tangent bundle $T\widetilde{M}$ splits into a direct
sum of holomorphic line bundles.

Since $T\widetilde{M}$ splits into a direct
sum of holomorphic line bundles, and $\widetilde{M}$ is
a compact connected K\"ahler manifold, using \cite[page 242,
Theorem 1.2]{BPZ} we know that $\widetilde{M}$ is biholomorphic
to the Cartesian product $({\mathbb C}{\mathbb P}^1)^d$.

Fix a point $x_0\, \in\, {\mathbb C}{\mathbb P}^1$.
Consider the map
$$
\phi\, :\, {\mathbb C}{\mathbb P}^1\, \longrightarrow\,
({\mathbb C}{\mathbb P}^1)^d\, =\, \widetilde{M}
$$
defined by $x\, \longmapsto\, (x\, , x_0\, ,\cdots\, ,x_0)$.
We have
$$
\phi^*T\widetilde{M}\, =\, {\mathcal O}_{{\mathbb C}{\mathbb
P}^1}(2) \bigoplus ({\mathcal O}_{{\mathbb C}{\mathbb
P}^1})^{\oplus (d-1)}\, .
$$
Since $d\, \geq\, 2$, it follows immediately from this
decomposition that the vector bundle $\phi^*T\widetilde{M}$ is
not semistable.

This contradicts the earlier observation that $\phi^*T
\widetilde{M}$ is semistable. Hence the fundamental group
of $M$ is infinite. This completes the proof of the theorem.
\end{proof}

\section{Maps from the projective line}

Let $M$ be a connected complex projective manifold
of complex dimension $d$, with $d\, \geq\, 2$. Assume that
$M$ has the property that for every pair of the
form $(Y\, ,\varphi)$,
where $Y$ is a compact connected Riemann surface and
$$
\varphi\,:\, Y\,\longrightarrow\, M
$$
a holomorphic map, the pull back $\varphi^*TM$ is a semistable
vector bundle over $Y$.

\begin{proposition}\label{prop2}
There is no nonconstant morphism from ${\mathbb C}{\mathbb P}^1$
to $M$.
\end{proposition}

\begin{proof}
To prove by contradiction, let
\begin{equation}\label{f}
f\, :\, {\mathbb C}{\mathbb P}^1\, \longrightarrow\, M
\end{equation}
be a nonconstant morphism. The given condition on $M$ says that
$f^*TM$ is a semistable vector bundle over ${\mathbb C}{\mathbb
P}^1$. Any holomorphic vector bundle over ${\mathbb C}{\mathbb
P}^1$ splits into a direct sum of holomorphic of line bundles
\cite[page 122, Th\'eor\`eme 1.1]{Gr}. Therefore, we have
\begin{equation}\label{f2}
f^*TM\, =\,
({\mathcal O}_{{\mathbb C}{\mathbb P}^1}(n))^{\oplus d}
\end{equation}
for some integer $n$.

The differential $df\, :\, T{\mathbb C}{\mathbb P}^1\, 
\longrightarrow\, f^*TM$ of the map $f$ in Eq. \eqref{f}
does not vanish identically because
$f$ is nonconstant. Since have a nonzero homomorphism from
$T{\mathbb C}{\mathbb P}^1\, =\, {\mathcal O}_{{\mathbb
C}{\mathbb P}^1}(2)$ to $f^*TM$, it follows that
$$
n\, \geq\, 2\, ,
$$
where $n$ is the integer in Eq. \eqref{f2}.
Consequently, the pull back $f^*TM$ is an ample vector bundle.

Since $f^*TM$ is ample, the variety $M$ is rationally connected
(see \cite[page 433, Theorem 2.1]{KMM} and \cite[page 434,
Definition--Remark 2.2]{KMM}). This in turn implies
that $M$ is simply connected \cite[p. 545, Theorem 3.5]{Ca},
\cite[page 362, Proposition 2.3]{Ko}. But this contradicts
Theorem \ref{th1}.

Therefore, there is no nonconstant morphism from
${\mathbb C}{\mathbb P}^1$ to $M$. This completes the proof
of the proposition.
\end{proof}

\section{The case of K\"ahler--Einstein manifolds}

As before, $M$ is a complex projective manifold of complex
dimension $d$, with $d\,\geq\, 2$.
Assume that there exists a K\"ahler form $\omega$ on $M$ with the 
following property:

There is a non--positive real number $\lambda\, \in\, {\mathbb C}$
such that the cohomology class of $\lambda\cdot\omega$ coincides
with the Chern class $c_1(TM)\, \in\, H^2(M,\, {\mathbb R})$.

A theorem due to Yau, \cite{Ya}, says that
there is a K\"ahler metric $\widetilde{\omega}$
on $M$ satisfying the following two conditions:
\begin{enumerate}
\item the K\"ahler metric $\widetilde{\omega}$
is K\"ahler--Einstein,
and \item the cohomology class $[\widetilde{\omega}]
\, \in\, H^2(M,\, {\mathbb R})$ coincides with that of $\omega$.
\end{enumerate}
(In \cite{Au}, this was proved under the
assumption that $c_1(TM)$ is positive.)

\begin{theorem}\label{th2}
Assume that for every compact connected Riemann surface $Y$, and
for every holomorphic map
$$
\varphi\, :\, Y\, \longrightarrow\, M\, ,
$$
the pull back $\varphi^*TM$ is a semistable vector bundle over $Y$.
Then $M$ admits a flat K\"ahler metric.
\end{theorem}

\begin{proof}
Let $\widetilde{\omega}$ be the K\"ahler--Einstein metric on
$M$. We will show that $\widetilde{\omega}$ is projectively flat.

Consider the Hermitian structure $\widetilde{\omega}'$
on the vector bundle $\mathcal{E}nd(TM)$ induced
by the Hermitian metric $\widetilde{\omega}$ on $TM$. Since
$\widetilde{\omega}$ is a K\"ahler--Einstein metric, it follows
that $\widetilde{\omega}'$ is a Hermitian--Einstein
metric. From Lemma \ref{lem1} we know that
$$
c_2(\mathcal{E}nd(TM))\, =\, 2d\cdot c_2(TM)-
(d-1)c_1(TM)^2 \, \in\, H^2(M, \, {\mathbb Q})
$$
vanishes. In view of this, the condition that $\widetilde{\omega}'$
is a Hermitian--Einstein metric implies that $\widetilde{\omega}'$
is flat (see \cite[Ch.~IV, page 115, Theorem 4.11]{Ko}).
Therefore, $\widetilde{\omega}$ is projectively flat.

Let
$$
\gamma\, :\, \widetilde{M}\, \longrightarrow\, M
$$
be a universal cover $M$. The pulled back K\"ahler metric
$\gamma^*\widetilde{\omega}$ on $\widetilde{M}$ is projectively
flat because $\widetilde{\omega}$ is projectively flat.
Consequently, the holonomy of $\gamma^*\widetilde{\omega}$
is contained in the center $U(1)\, \subset\, U(d)$, where
$d\, =\, \dim_{\mathbb C}M$.

Since the holonomy of $\gamma^*\widetilde{\omega}$ is contained
in the center $U(1)\, \subset\, U(d)$, and $\widetilde{M}$
is simply connected, we conclude the following. There are
connected Riemann surfaces $C_i$, $1\leq i\leq d$,
equipped with K\"ahler forms $\omega_i$, such that the
product K\"ahler manifold
$$
\prod_{i=1}^d (C_i\, , \omega_i)\, =\, (\prod_{i=1}^d C_i\, ,
\bigoplus_{i=1}^d \omega_i)
$$
is holomorphically
isometric to $\widetilde{M}$ equipped the K\"ahler form
$\gamma^*\widetilde{\omega}$ \cite[page 49, Theorem 3.2.7]{Jo}.

Fix a point $(x_1\, , \cdots \, ,x_{d-1})
\, \in\, \prod_{i=1}^{d-1} C_i$. Consider the holomorphic
Hermitian vector bundle of rank $d$ on $C_d$ obtained by
restricting to
$$
(x_1\, , \cdots \, ,x_{d-1}\, , C_d)\, \subset\,
\prod_{i=1}^{d} C_i
$$
the tangent bundle $T(\prod_{i=1}^d C_i)$ equipped with
the Hermitian metric $\bigoplus_{i=1}^d \omega_i$. We note that
this holomorphic Hermitian vector bundle on $C_d$ decomposes
into a direct sum of $(TC_d\, , \omega_d)$ with the trivial
holomorphic Hermitian vector bundle of rank $d-1$ with fiber
$\bigoplus_{i=1}^{d-1} T_{x_i} C_i$ equipped with the
Hermitian structure $\bigoplus_{i=1}^{d-1} \omega_i (x_i)$.

We observed earlier
that $\gamma^*\widetilde{\omega}$ is projectively flat.
Since the restriction of the holomorphic Hermitian vector
bundle $(T(\prod_{i=1}^d C_i)\, , \bigoplus_{i=1}^d \omega_i)$
to $C_d$ decomposes into the direct sum of $(TC_d\, , \omega_d)$
with a trivial holomorphic Hermitian vector bundle of
positive rank, it follows that $(TC_d\, , \omega_d)$ is a flat
line bundle. Indeed, the condition that the restriction of
$(T(\prod_{i=1}^d C_i)\, , \bigoplus_{i=1}^d \omega_i)$
to $C_d$ is projectively flat implies that curvatures on $C_d$
of all the direct summands coincide.

Interchanging $(TC_d\, , \omega_d)$ with $(TC_i\, ,
\omega_i)$, $i\, \in\, [1\, ,d-1]$, we conclude that
$$
\prod_{i=1}^d (C_i\, , \omega_i)
$$
is a flat K\"ahler manifold. Therefore, $\gamma^*
\widetilde{\omega}$
is a flat metric. This completes the proof of the theorem.
\end{proof}

\begin{corollary}\label{cor1}
Let $M$ be a connected complex projective manifold
of complex dimension $d$, with $d\, \geq\, 2$. Assume that
either the canonical line bundle $K_M$ is ample or $K^{-1}_M$
is ample or $c_1(K^{-1}_M)\, \in\, H^2(M,\, {\mathbb Q})$
vanishes. Then exactly one of the following two statements
is valid:
\begin{enumerate}
\item There is a pair $(Y\, ,\varphi)$,
where $Y$ is a compact connected Riemann surface and
$$
\varphi\,:\, Y\,\longrightarrow\, M
$$
a holomorphic map, such that the pull back $\varphi^*TM$ is
not semistable.

\item There is an \'etale covering map $A\, \longrightarrow\, M$,
where $A$ is an abelian variety.
\end{enumerate}
\end{corollary}

\begin{proof}
If $A\, \longrightarrow\, M$ is an \'etale covering,
where $A$ is an abelian variety, then the fiber
product $Y\times_M A$ is an \'etale
covering of $Y$, and furthermore, the pull back of
$\varphi^*TM$ to $Y\times_M A$ is a trivial vector bundle.
This implies that $\varphi^*TM$ is a semistable vector bundle.
Hence the two statements in the corollary cannot be
simultaneously valid.

If $K^{-1}_M$ is ample, then $M$ is rationally connected
\cite[p. 766, Theorem 0.1]{KMM2}. Hence in this case
Proposition \ref{prop2} implies that the first statement holds.

Now assume that either $K_M$ is
ample or $c_1(K_M)\, \in\, H^2(M,\, {\mathbb Q})$
vanishes. Then we know that
$M$ admits a flat K\"ahler metric (see Theorem \ref{th2}). 
If $M$ admits a flat K\"ahler metric, then the universal cover of
$M$ is ${\mathbb C}^d$, and the deck transformations are contained
in the automorphisms of ${\mathbb C}^d$ that preserve the constant
metric on ${\mathbb C}^d$. Consequently, the second statement in
the corollary holds. This completes the proof of the corollary.
\end{proof}

\begin{corollary}\label{cor2}
Let $M$ be a compact connected K\"ahler manifold such that the
rank of the N\'eron--Severi group
$\text{NS}(M)\, =\, H^2(M,\, {\mathbb Z})
\bigcap H^{1,1}(M)$ is one. Then exactly one of the following
two statements is valid:
\begin{enumerate}
\item There is a pair $(Y\, ,\varphi)$,
where $Y$ is a compact connected Riemann surface and
$$
\varphi\,:\, Y\,\longrightarrow\, M
$$
a holomorphic map, such that the pull back $\varphi^*TM$ is
not semistable.

\item There is an \'etale covering map $A\, \longrightarrow\, M$,
where $A$ is an abelian variety.
\end{enumerate}
\end{corollary}

\begin{proof}
If $\text{rank}(\text{NS}(M))\, =\, 1$, then
either $K_M$ is ample or $K^{-1}_M$ is ample
or $c_1(K_M)\, \in\, H^2(M,\, {\mathbb Q})$
vanishes. Therefore, the corollary follows from Corollary
\ref{cor1}.
\end{proof}

\section{The case of a surface}

Let $M$ be an irreducible smooth complex projective surface.

\begin{proposition}\label{prop3}
Exactly one of the following two statements is valid:
\begin{enumerate}
\item There is a pair $(Y\, ,\varphi)$,
where $Y$ is a compact connected Riemann surface and
$$
\varphi\,:\, Y\,\longrightarrow\, M
$$
a holomorphic map, such that the pull back $\varphi^*TM$ is
not semistable.

\item The surface $M$ admits an \'etale covering by an
abelian surface.
\end{enumerate}
\end{proposition}

\begin{proof}
It was shown in the proof of Corollary \ref{cor1} that
the two statements cannot be simultaneously valid.

Assume that the first statement does not hold. So
for every pair of the form $(Y\, ,\varphi)$,
where $Y$ is a compact connected Riemann surface and
$$
\varphi\,:\, Y\,\longrightarrow\, M
$$
a holomorphic map, the pull back $\varphi^*TM$ is a semistable
vector bundle over $Y$.

{}From Proposition \ref{prop2} we know
that $M$ is a minimal surface. If $M$ is of general type,
then $c_2(TM) \, >\, 0$, and also the Miyaoka inequality
$$
c_1(TM)^2 \, \leq\, 3c_2(TM)
$$
holds (see \cite[page 207, Theorem (1.1)]{BPV}). Hence
$$
c_1(TM)^2 - 4c_2(TM) \, <\, 0\, .
$$
This contradicts Lemma \ref{lem1}. Hence $M$ is not of
general type.

{}From Proposition \ref{prop2} we know that $M$ is not
a ruled surface.

Hence from the list of minimal projective surfaces (see
\cite[page 188, Table 10]{BPV}) we know $c_1(TM)^2 \, =\, 0$.
Therefore, from Lemma \ref{lem1} we know that
$c_2(TM) \, =\, 0$. Hence, from the list
of minimal projective surfaces we conclude that
exactly one of the following two statements holds:
\begin{enumerate}
\item The surface $M$ admits an \'etale covering by an
abelian surface.

\item There is an elliptic fibration $M\, \longrightarrow\, C$
with $\text{genus}(C) \, \geq\, 2$.
\end{enumerate}

The proof of the proposition will be completed
by showing that the second statement
does not hold. To prove this by contradiction, let
\begin{equation}\label{beta}
\beta\, :\, M\, \longrightarrow\, C
\end{equation}
be an elliptic fibration such that $C$ is a smooth projective
curve of genus at least two.

Since there is no nonconstant map from ${\mathbb C}{\mathbb P}^1$
to $M$, all the singular fibers of $\beta$ in Eq. \eqref{beta}
must be multiples of smooth elliptic curves. From this it
follows that there is a finite covering
$$
\alpha\, :\, \widetilde{C}\, \longrightarrow\, C
$$
with possible ramifications such that the normalization
$\widetilde{M}$ of the fiber product $M\times_C \widetilde{C}$
is a smooth elliptic fibration over $\widetilde{C}$, and
furthermore, the resulting morphism
\begin{equation}\label{b}
\gamma\, :\, \widetilde{M}\, \longrightarrow\, M
\end{equation}
is an \'etale covering map. Note that since
$$
c_1(TM)^2 \, =\, 0\, =\, c_2(TM)\, ,
$$
using the Hirzebruch--Riemann--Roch theorem it
follows that the Euler characteristic of
${\mathcal O}_M$ vanishes. Hence the above assertion
can be deduced using \cite[page 162, Remark]{BPV}
and \cite[page 162, Proposition (12.2)]{BPV}.

Since $\widetilde{M}\, \longrightarrow\, \widetilde{C}$
is a smooth elliptic fibration, the $j$--invariant map,
that associates to each point $x\, \in\, \widetilde{C}$
the $j$--invariant of the fiber $\widetilde{M}_x$ over $x$,
is in fact a constant map. Therefore, there is a finite
\'etale Galois covering
$$
\alpha'\, :\,\widetilde{C}'\,\longrightarrow\, \widetilde{C}
$$
such that
\begin{equation}\label{idf.}
\widetilde{M}'\, :=\, \widetilde{M}\times_{\widetilde{C}}
\widetilde{C}' \, =\, Z\times \widetilde{C}'\, ,
\end{equation}
where $Z$ is a smooth elliptic curve.

Let
\begin{equation}\label{bp}
\gamma'\, :\, \widetilde{M}'\, \longrightarrow\, M
\end{equation}
be the composition of the natural projection
$\widetilde{M}'\, \longrightarrow\, \widetilde{M}$
with the morphism $\gamma$ in Eq. \eqref{b}. Let
\begin{equation}\label{c}
\widehat{\gamma}\, :\, Z\times \widetilde{C}'\,
\longrightarrow\, M
\end{equation}
be the composition of the identification $\widetilde{M}'\,
=\, Z\times \widetilde{C}'$ in Eq. \eqref{idf.} with
$\gamma'$ in Eq. \eqref{bp}.

Fix a point $z_0\, \in\, Z$. Let
\begin{equation}\label{d}
\varphi_0\, :\, \widetilde{C}'\, \longrightarrow\,
\widetilde{M}' \, =\, Z\times \widetilde{C}'
\end{equation}
be the map defined by $c\,\longmapsto\, (z_0\, , c)$, where
$\widehat{\gamma}$ is constructed in Eq. \eqref{c}. Define
$$
\varphi\, :=\,\widehat{\gamma}\circ \varphi_0 
\, =\, \widehat{\gamma}\vert_{z_0\times \widetilde{C}'}
\, :\, \widetilde{C}'\, \longrightarrow\, M\, ,
$$
where $\widehat{\gamma}$ is constructed in Eq. \eqref{c}.
Since the map $\gamma'$ in Eq. \eqref{bp} is \'etale, we have
\begin{equation}\label{dc}
\varphi^*TM \, =\, \varphi^*_0 T\widetilde{M}'\, =\,
{\mathcal O}_{\widetilde{C}'}\bigoplus T\widetilde{C}'\, .
\end{equation}
Since $ \text{genus}(\widetilde{C}')\, \geq\,
\text{genus}(\widetilde{C})\, \geq\,
\text{genus}(C) \, \geq\, 2$, we have
$\text{degree}(T\widetilde{C}')\, \not=\, 0$. Hence from
Eq. \eqref{dc} it follows that $\varphi^*TM$ is not
semistable.

But this contradicts the initial assumption that
the first statement in the proposition does not hold.
Therefore, there is no elliptic fibration $M\, \longrightarrow\,
C$ with $\text{genus}(C) \, \geq\, 2$. This completes the proof
of the proposition.
\end{proof}


\end{document}